\documentclass[a4paper,12pt]{article}
\usepackage{latexsym,amsmath,amssymb,amsthm}
\theoremstyle{plain}
\newtheorem{theorem}{Theorem}[section]

\newtheorem{definition}[theorem]{Definition}

\date{}

\title{ End point of some generalized weakly contractive multivalued
mappings}
\author{Ali Abkar and Mohammad Eslamian}
\begin{document}
\maketitle Department of Mathematics, Imam Khomeini International
University, Qazvin 34149, Iran; email: abkar@ikiu.ac.ir,
mhmdeslamian@gmail.com

\begin{abstract}In this paper, we prove the existence of a common
end point for a pair of multivalued mappings satisfying a new
generalized weakly contractive condition in a complete metric space.
Our result generalizes and extends many known results.
\end{abstract}\maketitle\noindent {\bf Key words}: end point,
weakly contractive mapping, multivalued mapping.
\section{Introduction}
Banach contraction principle is a remarkable result in metric fixed
point theory. Over the years, tt has been generalized in different
directions and spaces by mathematicians.\par In 1997, Alber and
Guerre-Delabriere \cite{al} introduced the concept of weak
contraction in the following way:
\begin{definition} Let $(E,d)$ be a metric space. A mapping $T:E\to E$ is said to be
weakly contractive provided that $$d(Tx,Ty)\le
d(x,y)-\varphi(d(x,y))$$
 where $x,y\in E$ and $\varphi:[0,\infty)\to\:[0,\infty)$ is  a
continuous nondecreasing function such that $\varphi(t)=0$ if and
only if $t=0$.
\end{definition}
Using the concept of weakly contractiveness, they succeeded to
establish the existence of fixed points for such mappings in Hilbert
spaces. Later on Rhoades \cite{a9} proved that the result of
\cite{al} is also valid in complete metric spaces. Rhoades \cite{a9}
also proved the following fixed point theorem which is a
generalization of the Banach contraction principle, because it
contains contractions as a special case when we assume that
$\varphi(t)=(1-k)t$ for some $0< k <1$.
\begin{theorem} Let $(E,d)$ be a complete metric space and let $T:E\to E$ be a weakly
contractive mapping. Then $T$ has a fixed point.
\end{theorem}
Since then, fixed point theory for single valued, as well as for
multivalued weakly contractive type mappings was studied by many
authors; see [2-8], and [10-12]. \par Let $(E,d)$ be a metric space,
and let $B(E)$ denote the family of all nonempty bounded subsets of
$E$. For $A,B\in B(E)$, define the distance between $A$ and $B$ by
$$D(A,B)=\inf \{d(a,b): a\in A,b\in B\},$$ and the diameter of $A$
and $B$ by
$$\delta(A,B)=\sup\{ d(a,b): a\in A,b\in B\}.$$ Let $T:E\to B(E)$
be a multivalued operator, then an element $x\in E$ is said to be a
\emph{fixed point} of $T$ provided that $x\in T(x)$ and it is called
an \emph{end point} of $T$ if $T(x)=\{x\}$. The purpose of this
paper is to prove the existence of a common end point for a pair of
multivalued mappings satisfying a new generalized weakly contractive
condition in a complete metric space. Our result generalizes and
extends some a number of known results.
\section{Main Results}
In this sequel, we denote by $\Phi$ the class of all mappings
$\varphi:[0,\infty)\to [0,\infty)$ satisfying the following
conditions:\begin{enumerate}
\item [(i)] $\varphi (t)=0$ if and only if $t=0$;
\item [(ii)] $\varphi$ is a lower semi continuous function ;
\item [(iii)] for any sequence $ \{t_{n}\}$ with
$\lim_{n\to\infty}t_{n}=0$ , there exist $k\in(0,1)$ and $n_{0}\in
\mathbb{N}$, such that $\varphi(t_{n})\geq kt_{n}$ for each $n\geq
n_{0}.$\end{enumerate} Examples of such mappings are $\varphi
(x)=kx$ for $0<k<1$ and $\varphi(x)=\ln(x+1)$ (see also \cite{a8}).
Let $\Omega$ denote the class of all mappings $f:[0,\infty)\to
[0,\infty)$ satisfying the following conditions:\begin{enumerate}
\item [(i)] $f(t)=0$ if and only if $t=0$;
\item [(ii)] $f$ is  non-decreasing;
\item [(iii)] $f$ is continuous;
 \item [(iv)]  $f(x+y)\leq f(x)+f(y)$.\end{enumerate}
Finally, let $\Psi$ denote the class of mappings $\psi:[0,\infty)\to
[0,\infty)$ which are continuous and non-decreasing with $\psi(t)=0$
if and only if $t=0$.\par
 Let $(E,d)$ be a metric space, and let $T,S:E\to B(E)$
 be two multivalued  mappings, we define
$$ M(x,y)=\max\left
\{d(x,y),\delta(Tx,x),\delta(y,Sy),\frac{D(y,Tx)+D(x,Sy)}{2}\right
\},$$ and $$N(x,y)=\min\{D(y,Tx),D(x,Sy)\}.$$ We now state the main
result of this paper.
\begin{theorem}Let $(E,d)$ be a complete metric space, and let $T,S:E\to B(E)$
 be two mappings such that for all $x,y\in E$
\begin{equation}
f(\delta(Tx,Sy)\le f(M(x,y))-\varphi(f(M(x,y)))+\psi(N(x,y))
\end{equation}
where $\varphi\in \Phi$, $\psi\in\Psi$ and $f\in \Omega.$ Then $S$
and $T$ have a common end point $z \in E$, i.e, $Sz=Tz=\{z\}$.
\end{theorem}
\begin{proof}We construct a sequence $\{x_{n}\}$ as follows. Take $x_{0}\in E$ and for $n\geq 1$ we choose
$ x_{2n+1}\in Tx_{2n}:=A_{2n}$ and  $x_{2n+2}\in
Sx_{2n+1}:=A_{2n+1}.$ Now we have
\begin{multline*}M(x_{2n},x_{2n+1})\\= \max
\{d(x_{2n},x_{2n+1}),\delta(Tx_{2n},x_{2n}),\delta(Sx_{2n+1},x_{2n+1}),\\
\frac{D(Tx_{2n},x_{2n+1})+D(Sx_{2n+1},x_{2n})}{2}\}\\
\leq \max\{\delta( A_{2n-1},A_{2n}),\delta(
A_{2n-1},A_{2n}),\delta(A_{2n+1},A_{2n}),\\
\frac{D(Tx_{2n},x_{2n+1})+D(Sx_{2n+1},x_{2n})}{2} \}\\\leq\max
\{\delta(
A_{2n-1},A_{2n}),\delta(A_{2n+1},A_{2n}),\frac{\delta(A_{2n+1},A_{2n-1})}{2}\}\\\leq
\max\{\delta(A_{2n-1},A_{2n}),\delta(A_{2n+1},A_{2n}),\\
\frac{\delta(A_{2n},A_{2n-1})+\delta(A_{2n},A_{2n+1})}{2}\}\\
=\max\{\delta(A_{2n-1},A_{2n}),\delta(A_{2n+1},A_{2n})\}
\end{multline*}
and $$N(x_{2n},x_{2n+1})\\= \min
\{D(Tx_{2n},x_{2n+1}),D(Sx_{2n+1},x_{2n}) \}=0.$$ By assumption
\begin{multline*}
f(\delta (A_{2n},A_{2n+1}))=f(\delta (Tx_{2n},Sx_{2n+1}))\\
\leq
f(M(x_{2n},x_{2n+1}))-\varphi(f(M(x_{2n},x_{2n+1})))+\psi(N(x_{2n},x_{2n+1}))
\\=f(M(x_{2n},x_{2n+1}))-\varphi(f(M(x_{2n},x_{2n+1}))\\\leq
f(M(x_{2n},x_{2n+1}).
\end{multline*}
Since $f$ is non-decreasing, we have
$$\delta (A_{2n},A_{2n+1})\leq M(x_{2n},x_{2n+1}).$$
Now, if $\delta(A_{2n-1},A_{2n})<\delta(A_{2n+1},A_{2n})$ then
$$M(x_{2n},x_{2n+1})\leq \delta (A_{2n+1},A_{2n}),$$
from which we obtain $$M(x_{2n},x_{2n+1})= \delta
(A_{2n+1},A_{2n})>\delta(A_{2n-1},A_{2n})\geq 0,$$ and
\begin{multline*}
f(\delta (A_{2n},A_{2n+1}))=f(\delta (Tx_{2n},Sx_{2n+1}))\\
\leq
f(M(x_{2n},x_{2n+1}))-\varphi(f(M(x_{2n},x_{2n+1})))+\psi(N(x_{2n},x_{2n+1}))
\\=f(M(x_{2n},x_{2n+1}))-\varphi(f(M(x_{2n},x_{2n+1}))\\<
f(M(x_{2n},x_{2n+1})=f(\delta(A_{2n+1},A_{2n}))
\end{multline*}
which is a contradiction. So we have
$$\delta(A_{2n+1},A_{2n})\leq M(x_{2n},x_{2n+1})\leq
\delta(A_{2n},A_{2n-1}).$$
 Similarly we obtain
$$\delta(A_{2n+1},A_{2n+2})\leq M(x_{2n+1},x_{2n+2})\leq
\delta(A_{2n+1},A_{2n+2}).$$ Therefore the sequence
$\{\delta(A_{n},A_{n+1})\}$ is monotone decreasing and bounded
below. So there exists $r\geq 0$ such that
$$\lim_{n\to\infty}{\delta(A_{n},A_{n+1})}=\lim_{n\to\infty}{M(x_{n},x_{n+1})}=r.$$
We now claim that $r=0$. In fact taking upper limits as
$n\to\infty$ on either sides of the inequality
\begin{multline*}
f(\delta (A_{2n},A_{2n+1}))=f(\delta (Tx_{2n},Sx_{2n+1}))\\
\leq
f(M(x_{2n},x_{2n+1}))-\varphi(f(M(x_{2n},x_{2n+1})))+\psi(N(x_{2n},x_{2n+1}))
\\=f(M(x_{2n},x_{2n+1}))-\varphi(f(M(x_{2n},x_{2n+1})))
\end{multline*}
we have $$f(r)\leq f(r)-\varphi f(r)$$ which is a contradiction
unless $r=0$. Thus $\lim_ {n\to\infty}\delta(A_{n},A_{n+1})=0$
 and hence $\lim_ {n\to\infty}d(x_{n},x_{n+1})=0.$
Now we shall prove that $\{x_n\}$ is a Cauchy sequence. Indeed,
Since $\lim_ {n\to\infty}f(M(x_{n},x_{n+1}))=0$ by the property of
$\varphi$ there exist $0<k<1$ and $n_{0}\in\mathbb{N}$, such that
$\varphi(f(M(x_{n},x_{n+1})))\geq k f(M(x_{n},x_{n+1}))$ for all
$n>n_{0}.$ On the other hand, for any given $\varepsilon>0$, we can
choose $\eta>0$ in such a way that $f(\eta)\leq\frac
{k}{1-k}f(\varepsilon).$ Moreover, there exists $n_{0}$ such that
$\delta(A_{n},A_{n-1})\leq\eta$ for each $n>n_{0}$. For any natural
number $m>n>n_{0}$ if $n$ is even, we have
\begin{multline*}
f(\delta(A_{n},A_{n+1}))\leq f(\delta(Tx_{n},Sx_{n+1}))\\\leq
f(M(x_{n},x_{n+1}))-\varphi(f(
M(x_{n},x_{n+1}))+\psi(N(x_{n},x_{n+1}))\\\leq (1-k)f(
M(x_{n},x_{n+1}))\leq (1-k) f(\delta(A_{n},A_{n-1})).
\end{multline*}
By this inequality, we get for $l>n$
$$f(\delta(A_{l},A_{l-1}))\leq (1-k)
f(\delta(A_{l-1},A_{l-2}))\leq \cdots \leq(1-k)^{l-n}
f(\delta(A_{n},A_{n-1}))$$ Therefore we have
\begin{multline*}f(\delta(A_{n},A_{m}))\leq
f(\delta(A_{n},A_{n+1})+\delta(A_{n+1},A_{n+2})+\cdots+\delta(A_{m-1},A_{m}))\\
\leq
f(\delta(A_{n},A_{n+1}))+f(\delta(A_{n+1},A_{n+2}))+\cdots+f(\delta(A_{m-1},A_{m}))\\
\leq (1-k)f(\delta(A_{n},A_{n-1}))+\cdots\\
+(1-k)^{m-n-1}f(\delta(A_{n},A_{n-1}))+(1-k)^{m-n}f(\delta(A_{n},A_{n-1}))\\
=\frac{(1-k)-(1-k)^{m-n+1}}{1-(1-k)}f(\delta(A_{n},A_{n-1}))\\
< \frac{1-k}{k}f(\delta ((A_{n},A_{n-1})))\leq
\frac{1-k}{k}f(\eta)<f(\varepsilon) .
\end{multline*} Now, by the nondecreasingness of $f$ we obtain $\delta(A_{n},A_{m})<
\varepsilon$. From the construction of the sequence $\{x_{n}\}$, it
follows that the same conclusion holds for $\{x_{n}\}$, i.e. for
each $\varepsilon>0$ there exist $n_{0}$ such that for any natural
numbers $m>n>n_{0}$, $ d(x_{n},x_{m})<\varepsilon$. This shows that
$\{x_{n}\}$ is a Cauchy sequence. Notice that $E$ is complete, hence
$\{x_n\}$ is convergent. Let us denote its limit by
$\lim_{n\to\infty} x_{n}=z$ for some $z\in E$. Now we prove that
$\delta(Tz,z)=0$. Suppose that this is not true, then
$\delta(Tz,z)>0$. For large enough $n$, we claim that the following
equations hold true:
\begin{multline*}M(z,x_{2n+1})
=\max \{d(z,x_{2n+1}),\delta(z,Tz),
\delta(Sx_{2n+1},x_{2n+1}),\\
\frac{D(Tz,x_{2n+1})+D(Sx_{2n+1},z)}{2}\}
 =\delta(z,Tz).\end{multline*}
 Indeed, since  $$\delta(Sx_{2n+1},x_{2n+1})\leq \delta(A_{2n+1},A_{2n})\to
 0,$$
 and
\begin{multline*}\lim_{n\to\infty}{\frac{D(Tz,x_{2n+1})+D(Sx_{2n+1},z)}{2}}\\
\le\lim_{n\to\infty}{\frac{\delta(Tz,z)+d(z,x_{2n+1})+\delta(Sx_{2n+1},x_{2n+1})+d(x_{2n+1},z)}{2}}\\
=\frac{\delta(Tz,z)}{2},
\end{multline*}
it follows that there exists $k\in \mathbb{N}$ such that
$M(z,x_{2n+1})=\delta(z,Tz)$ for $n>k$. Note that
\begin{multline*}f(\delta(Tz,x_{2n+2}))\le
f(\delta(Tz,Sx_{2n+1}))\\
\le f(M(z,x_{2n+1}))-\varphi(f((M(z,x_{2n+1}))-\psi(N(z,x_{2n+1}))
.\end{multline*}
 Letting $n\to\infty$, we have
$$f(\delta(Tz,z))\le f(\delta(Tz,z))-\varphi(f(\delta(Tz,z)))$$
i.e, $\varphi(f(\delta(Tz,z)))\le 0$. This is a contradiction,
therefore $\delta(Tz,z)=0$ i.e., $Tz=\{z\}$. And
since\begin{multline*}M(z,z)=\max\{d(z,z),\delta(Tz,z),
\delta(z,Sz),\frac{D(Tz,z)+D(Sz,z)}{2}\}\\
=\max\{\delta(Sz,z),\frac{D(Sz,z)}{2}\}=\delta(Sz,z)\end{multline*}
and $$N(z,z)=\min\{D(z,Tz),D(z,Sz)\}=0$$
we conclude that \begin{multline*}f(\delta(z,Sz))\le f(\delta(Tz,Sz))\\
\le f(M(z,z))-\varphi(f(M(z,z)))+\psi (N(z,z))\\
\le f(\delta(z,Sz))-\varphi(f(\delta(Sz,z))),\end{multline*}which in
turn implies that $Sz=\{z\}$. Hence the point $z$ is a common end
point of $S$ and $T$.\\
\end{proof}
\begin{theorem}Let $(E,d)$ be a complete metric space, and let $T,S:E\to B(E)$
 be two mappings such that for all $x,y\in E$
\begin{equation}
f(\delta(Tx,Sy)\le f(M(x,y))-\varphi(f(M(x,y)))
\end{equation}
where $\phi\in \Phi$ and $f\in \Omega.$ Then $S$ and $T$ have a
unique common end point $z \in E$.i.e, $Sz=Tz=\{z\}$.
\end{theorem}
\begin{proof}
By theorem 2.1, $T$  and $S$ have a common end point $z$. Now let
$y\in E$ be another common end point of $S$ and $T$. Notice
that\begin{multline*}M(y,y)=\max\{d(y,y),\delta(Ty,y),\delta(y,Sy),
\frac{D(Ty,y)+D(Sy,y)}{2}\}\\
=\max\{\delta(Sy,y),\delta(y,Ty)\}.\end{multline*}
Hence\begin{multline*}f(\delta(y,Ty))\le f(\delta(Sy,Ty))
\le f(M(y,y))-\varphi(f(M(y,y)))\\
\le f
(max\{\delta(y,Sy),\delta(y,Ty)\})-\varphi(f(\max\{\delta(y,Sy),\delta(y,Ty)\})).\end{multline*}
Similarly, we have
\begin{multline*}f(\delta(y,Sy))\le f(\delta(Ty,Sy))
\le f(M(y,y))-\varphi(f(M(y,y)))\\
\le f
(\max\{\delta(y,Sy),\delta(y,Ty)\}-\varphi(f(\max(\delta(y,Sy),\delta(y,Ty)\})).\end{multline*}
Therefore \begin{multline*}f(\max\{\delta(y,Sy),\delta(y,Ty)\})\le
\\
f(\max\{\delta(y,Sy), \delta(y,Ty)\})-\phi
(f(\max\{\delta(y,Sy),\delta(y,Ty)\})))\end{multline*} which implies
that $\max\{\delta(y,Sy),\delta(y,Ty)\}=0$, hence
$\delta(Ty,y)=\delta(Sy,y)=0.$ Now we
have$$M(z,y)=\max\{d(z,y),\delta(z,Tz),\delta(y,Sy),\frac{D(y,Tz)+D(z,Sy)}{2}\}$$
and
\begin{multline*}f(d(z,y))=f(\delta(Sz,Ty))\le f(M(z,y))-\varphi(f(M(z,y)))\\
=f(d(z,y))-\varphi (f(d(z,y)))\end{multline*} that imply
$d(z,y)=0$ i.e, $z=y$. Hence $z$ is the unique common end point of
$S$ and $T$.
\end{proof}
If in Theorem 2.1 we put $f(t)=t$ and $\varphi(t)=(1-k)t$, for some
$0<k<1$, then we obtain the following result.
\begin{theorem}Let $(E,d)$ be a complete metric space, and let $T,S:E\to B(E)$
 be two mappings such that for all $x,y\in E$
\begin{equation}
\delta(Tx,Sy)\le k(M(x,y))+\psi(N(x,y))
\end{equation}
where $\psi\in\Psi$. Then $S$ and $T$ have a common end point $z \in
E$, i.e, $Sz=Tz=\{z\}$.
\end{theorem} Let $T$ and $S$ be two single valued mappings, the we obtain the following theorem:
\begin{theorem}Let $(E,d)$ be a complete metric space, and let $T,S:E\to E$
 be two mappings such that for all $x,y\in E$
\begin{equation}
f(d(Tx,Sy)\le f(M(x,y))-\varphi(f(M(x,y)))+\psi(N(x,y))
\end{equation}
where $\varphi\in \Phi$, $\psi\in\Psi$, \, $f\in \Omega$ and
$$M(x,y)=\max\left
\{d(x,y),d(Tx,x),d(y,Sy),\frac{d(y,Tx)+d(x,Sy)}{2}\right \},$$
$$N(x,y)=\min \{d(Tx,y),d(x,Sy)\}.$$ Then $S$ and $T$ have a common
fixed point $z \in E$, i.e, $Sz=Tz=z$.
\end{theorem}

\noindent{\bf Example 1}. Let $E=[0,1]$ and $d(x,y)=|x-y|$. For each
$x\in E$ define $S,T:E\to B(E)$ by
$$Tx=[\frac{x}{4},\frac{x}{2}] ,\qquad Sx=[0,\frac{x}{5}].$$ Then
$$\delta(Tx,Sy)=\begin{cases}\frac{x}{2} &0\le\frac {y}{5}\le
\frac{x}{2}\\ \max\{\frac{y}{5}-\frac{x}{4},\frac{x}{2}\} &\frac x2
\le\frac y5\le 1.\end{cases}$$ and
$$\delta(x,Tx)=\frac{3x}{4},\qquad \delta(y,Sy)=y.$$ We also consider
$f(t)=2t$ and $\phi(t)=\frac{t}{4}$. We note that if
$\delta(Tx,Sy)=\frac{x}{2}$
then\begin{multline*}f(\delta(Tx,Sy)=x\le\frac{9x}{8}=\frac{3}{2}\delta(x,Tx)\\
\le\frac{3}{2}(M(x,y))=f(M(x,y))-\varphi(f(M(x,y)))\end{multline*}
 and if $\delta(Tx,Sy)=\frac{y}{5}-\frac{x}{4}$ then\begin{multline*}
\psi(\delta(Tx,Sy)=2(\frac{y}{5}-\frac{x}{4})\le\frac{2y}{5}\le\frac{3y}{2}\\
=\frac{3}{2}\delta(y,Sy)\le\frac{3}{2}(M(x,y))
=f(M(x,y))-\varphi(f(M(x,y))). \end{multline*} This arguments show
that the mappings $T$ and $S$ satisfy the conditions of Theorem 2.2.
Now it is easy to see that $0$ is the only common end point of this
two mappings.\par In the following we shall see that Theorem 2.1 is
a real generalization of Theorem 2.2. We note that by Theorem 2.2,
$T$ and $S$ have a unique common end point.

\noindent{\bf Example 2}. Let $E=[0,1]$ and $d(x,y)=|x-y|$. For each
$x\in E$ define $T,S:E\to B(E)$ by
$$Tx=Sx=\begin{cases}[\frac{x}{3},\frac{x}{2}] & x\neq 1\\ 1  & 1.\end{cases}$$ We also consider
$f(t)=t$, $\phi(t)=\frac{t}{5}$, and $\psi(t)=2t^{2}.$
 We note that if $x,y\neq 1$ and $x\leq y$
then $\delta (y,Ty)=\frac{2y}{3}$ hence
 \begin{multline*}f(\delta(Tx,Ty)=\frac{y}{2}-\frac{x}{3}\leq \frac{8}{15}\delta(y,Ty)\\
\le\frac{8}{15}M(x,y)=f(M(x,y))-\varphi(f(M(x,y)))\end{multline*}
Similar result holds if $x,y\neq 1$ and $y\leq x.$ Now if  $y=1$ and
$x\neq1$, then  $\delta(Tx,Ty)=1-\frac{x}{3}$ ,
$D(y,Tx)=1-\frac{x}{2}$ , $D(x,Ty)=1-x$ , $d(x,y)=1-x$,
$\delta(y,Ty)=0$ and $\delta(x,Tx)=\frac{2x}{3}.$ Hence
$$M(x,y))=\begin{cases}1-\frac{3x}{4}, & x\leq\frac{12}{17} \\
\frac{2x}{3},  & x\geq\frac{12}{17} .\end{cases}$$  and
$$N(x,y)=1-x$$ therefore we have
\begin{multline*}f(\delta(Tx,Ty))=1-\frac{x}{3}\leq\\
2(1-x)^{2}=\psi(N(x,y))\leq
f(M(x,y))-\varphi(f(M(x,y)))+\psi(N(x,y)).\end{multline*}Similar
result holds if $x=1$ and $y\neq1$ This arguments show that the
mappings $T$ and $S$ satisfy the conditions of Theorem 2.1. We
observe that $0$ and $1$ are two end points for $T$ and $S$.

\end{document}